\documentclass[reqno, oneside]{amsart}
\usepackage{amsmath, amscd, amssymb}

\theoremstyle{plain}
\newtheorem{Thm}{Theorem}
\newtheorem{Prop}[Thm]{Proposition}
\newtheorem{Lem}[Thm]{Lemma}
\newtheorem{Cor}[Thm]{Corollary}

\theoremstyle{definition}

\theoremstyle{remark}
\newtheorem*{Rem}{Remark}

\numberwithin{Thm}{section}

\def\fl(#1){\eqref{fl:#1}}
\def\thm#1{Theorem~\ref{thm:#1}}
\def\lem#1{Lemma~\ref{lem:#1}}
\def\cor#1{Corollary~\ref{cor:#1}}
\def\prop#1{Proposition~\ref{prop:#1}}
\def\sec#1{\S{\rm \ref{sec:#1}}}

\newcommand{\R}{\bar{\mathcal R}}
\newcommand{\Lg}{\frak g}

\newcommand{\Ug}{U_{v}(\frak g)}

\newcommand{\Uq}{U_{q}^{res}(\frak g)}

\newcommand{\Uf}{{U_q^{fin}}(\frak g)}
\newcommand{\bUf}{{\stackrel{\bullet}{U_q^{fin}}}(\frak g)}

\newcommand{\Ud}{U_q^{fin}(\frak g)^*}

\newcommand{\Z}{\frak Z}

\newcommand{\Wl}{\hat{W}_{l,Q}}
\newcommand{\WlP}{\hat{W}_{l,P}}

\newcommand{\Cf}{{\mathcal C}_f}
\newcommand{\bCf}{\bar{{\mathcal C}_f}}

\newcommand{\tmX}{\bar{\mathcal X}}
\newcommand{\mR}{\mathcal R}

\newcommand{\mX}{\mathcal X}
\newcommand{\bJ}{\bar{J}}
\newcommand{\bVr}{\bar{Vr}}
\newcommand{\bPr}{\bar{Pr}}

\begin{document}
\author{Anna Lachowska}
\thanks{This research was conducted by the author for the Clay Mathematics
 Institute}

\address{Department of Mathematics\\
MIT\\
77 Massachusetts ave. \\
Cambridge MA 02139\\
USA}
\email{lachowska@math.mit.edu}

\title[A Counterpart of the Verlinde Algebra]{A Counterpart of the 
Verlinde Algebra for the Small Quantum Group}

\date\today

\begin{abstract} 
Let $\bPr$ denote the ideal spanned by the characters of projective 
modules in the Grothendieck ring of  the category $\bCf$ of finite 
dimensional modules over the small quantum group $\Uf$. We show 
that $\bPr$ admits a description completely parallel to that of the 
Verlinde algebra of the fusion category \cite{AP}, with the character of the Steinberg module playing the role of the identity.

\end{abstract}
\maketitle

\section{Introduction and Notations} \label{sec:In}

For a semisimple complex Lie algebra $\Lg$ with the root system $R$ and the 
weight lattice $P$, consider the category $\mathcal O$ 
introduced in \cite{BGG}. For any $\nu \in P$ let $L(\nu)$ and $M(\nu)$ 
denote respectively the simple and 
the Verma
module of highest weight $\nu - \rho$, where $\rho$ is half the 
sum of all positive roots, and let $P(\nu)$ be 
the projective cover of $L(\nu)$ in $\mathcal O$. 
It is well known (e.g. \cite{Hum1}) 
that for a dominant integral $\lambda \in P_+$ and the 
longest element 
$w_0$ of the Weyl group $W$, the projective module $P(w_0(\lambda))$ 
has a particularly nice structure. It is usually called the big projective 
module and 
 has a Verma flag with the multiplicity of each 
$M(w(\lambda))$, $w \in W$, equal to one. 
In terms of characters, if we write  
$ch M = \sum_{\eta \in P}(dim M_\eta)e^\eta \in {\mathbb C}[P]$ for the 
formal 
character of $M \in {\mathcal O}$, then for a big projective module we have 
$$ ch P(w_0(\lambda)) = \sum_{w \in W/W_\lambda}ch M(w(\lambda)) = ch M(0) 
\frac{1}{|W_\lambda|} \sum_{w \in W} e^{w (\lambda)}, $$
where $W_\lambda \subset W$ is the stabilizer subgroup of the weight 
$\lambda$.  
The elements 
$$ f(\lambda) = \frac{1}{|W_\lambda|} \sum_{w \in W} e^{w (\lambda)}, 
\;\;\;\;\; \lambda \in P_+, $$ 
can be viewed as a natural basis in the space ${\mathbb C}[P]^W$ 
of exponential 
$W$-invariants on the weight lattice $P$. The relation between this basis and 
the basis of Weyl characters
$$ \chi(\lambda) = \frac{\sum_{w \in W} \varepsilon(w)e^{w(\lambda +\rho)}}
{\sum_{w \in W} \varepsilon(w)e^{w(\rho)}} = 
ch M(0) \sum_{w \in W} \varepsilon(w)e^{w(\lambda +\rho)},
 \;\;\;\;\; \lambda \in P_+, $$ 
reflects a reciprocity between finite dimensional simple and big 
projective modules in the category $\mathcal O$. 

We will investigate a similar reciprocity in the framework 
of the representation theory of quantum groups at roots 
of unity. More precisely, we will describe a certain finite set of 
projective modules over the 
small quantum group corresponding to $\Lg$, whose characters are 
given by the same formula as above with 
the character of the Steinberg module replacing $chM(0)$ 
(cf. \cite{Hum2}, \S 9).   
The tensor product induces a multiplication in the set 
of characters of these projective modules, which leads 
to a structure similar to the well known 
Verlinde algebra of characters 
of a certain set of simple modules over the quantum group at a root of unity. 
The comparative  study of these two algebraic objects is the main 
goal of the present paper. We briefly outline the results below; the precise 
definitions can be found at the end of this section.  

Let $Q = {\mathbb Z}R$ denote the root lattice of $\Lg$. We will 
assume $R$ to be irreducible and simply laced, 
although the latter restriction is not essential. 
Fix an integer $l \geq h$, where $h$ is the Coxeter number of $R$ 
(later we will require $l$ to be odd and coprime to the determinant of the 
Cartan matrix of $R$, \sec{fus}), and consider the 
shifted action of the affine Weyl group  $\Wl \simeq  W\ltimes lQ$ 
on $P$. Define the elements $\chi(\lambda)$ for nondominant weights 
by setting $\chi(\lambda) = \varepsilon(w) \chi(w \cdot \lambda)$ for any 
$\lambda \in P, w \in W$. Then the quotient 
$$ (Vr)^- \simeq {\mathbb C}[P]^W/\langle (\chi(\lambda) + \chi(s \cdot \lambda)\rangle_{\lambda \in P, s - {\rm a}\,{\rm reflection}\, {\rm in}\, \Wl},$$ 
is isomorphic to the well known Verlinde algebra $Vr$ (\cor{J}). On the other 
hand, one can consider the natural action of the extended affine group 
 $\WlP \simeq 
W\ltimes lP$, where the translations are generated by the 
$l$-multiple weight lattice. Set $f(w(\lambda))=f(\lambda)$ for any $\lambda \in P$, $w \in W$, and consider the quotient  
$$ (Vr)^+ \simeq {\mathbb C}[P]^W/\langle (f(\lambda) - f(s (\lambda))\rangle_{\lambda \in P, 
s - {\rm a}\,{\rm reflection}\, {\rm in}\, \WlP}$$
with respect to the natural action.  
The algebra $(Vr)^-$ has a basis of Weyl characters with highest weights 
parametrized by $X$, the regular elements of the fundamental 
domain of $\Wl \cdot$ action, while a basis in  $(Vr)^+$ is given by 
$f(\lambda)$ 
with $\lambda \in \hat{\mX}$, the fundamental domain of the natural action of 
$\WlP$. 
 The analogy between the algebraic structures and the roles played by the two bases are self-evident. We explain below their representation-theoretical meaning. 

 Originally  $Vr$ was introduced in \cite{Ver} as the fusion ring 
with a finite number of simple objects, associated to a semisimple 
Lie algebra $\Lg$ and a number $l \geq h$. The fusion category can be realized in many ways, 
in particular for our purposes we will view it as a subquotient of the category of finite dimensional modules over the restricted specialization of a quantum group at the $l$-th root of unity \cite{And}, \cite{AP}.

It turns out that $(Vr)^+$ is not only a close analog of $(Vr)^-$ from the purely algebraic point of view, but it also shares many important properties of the Verlinde algebra and 
has a clear representation-theoretical meaning. Namely, let $St$ denote the Steinberg module, which is a simple, projective and injective object in the category of finite dimensional modules over the big quantum group $\Uq$. 
Its character $ch St$ is given by the Weyl formula of the highest weight $(l-1)\rho$. We show (\cor{bPr}) that the structure of $(Vr)^+$ defines 
multiplication   in the ideal $\bPr$ of elements corresponding to projective 
modules in the Grothendieck ring $\R$ of finite dimensional representations 
over the small quantum group $\Uf$ in the following sense: 
$$ [P(\lambda)][P(\mu)] = [St] \sum_{\nu}n_{\lambda, \mu}^\nu [P(\nu)]. $$ 
Here elements $\{[P(\lambda)]\}_{\lambda \in \hat{\mX}} \in \R$ 
correspond to indecomposable projective modules, and the numbers 
$n_{\lambda, \mu}^\nu$ are the structure constants of the algebra 
$(Vr)^+$ in the basis of $W$-symmetric functions $\{f(\lambda)\}_{\lambda \in 
\hat{\mX}}$.    
We will also see that for odd $l >h$, coprime to the determinant of the Cartan matrix, the ideal $\bPr$  carries a projective representation of the modular 
group (\prop{fourier}), whereas the algebra $Vr$ has this property for even 
$l$.
  
While the Verlinde algebra has been extensively studied in the literature, the existence 
of its analog $\bPr$ seems to be new. Based on the analogy between $Vr$ 
and $\bPr$, we expect that the latter will have important applications in 
various fields where the Verlinde algebra is involved. 
In particular, the next natural step would be to investigate the possibility of constructing a tensor category with primary objects in one-to-one correspondence with the basis in $\bPr$. Another open question is to identify a tensor 
category over the affine Lie algebra $\hat{\Lg}$ which corresponds to the 
multiplicative structure of $\bPr$, similarly to the equivalence of the fusion 
categories stated in \cite{Fin}. The contravariant equivalence between the category with finite Weyl filtrations over $\Uq$ and an analogous category over $\hat{\Lg}$ at a fixed positive level 
\cite{Soe1}, provides a tool to identify the counterparts of the projective modules among the positive level representations of $\hat{\Lg}$, though it does not supply the tensor structure.

The paper is organized as follows: \sec{fus} contains the description of the 
Verlinde algebra \cite{AP}, and its restriction $\bVr$ over the small quantum 
group (\thm{bJ}).In \sec{Pr} we study the algebraic structures behind the category of projective modules over $\Uq$ (\thm{Pr}) and $\Uf$ (\thm{bPr}), and consider the obtained results in the framework of the classification of 
tensor ideals in the category of tilting modules $\mathcal T$ over 
$\Uq$ \cite{Ost}.  
Finally, in \sec{VrPr} we list the similarities of $Vr$ and $\bPr$.

\subsection*{Definitions and notations.}
We will assume that $l$ is odd everywhere except in the description of 
the modular group action on the Verlinde algebra at the end of \sec{VrPr}.
Let $\Lg$ be a simply laced simple complex algebra with the root system 
$R$ and the Weyl group $W$.
 The action of $W$ on  
$P$ is given by $s_\alpha(\mu) = \mu - \langle \mu, \alpha \rangle \alpha$ for any 
$\alpha \in R$, $\mu \in P$, where the inner product 
$\langle \cdot, \cdot \rangle$  is normalized so that $\langle \alpha, \alpha 
\rangle =2$ for $\alpha \in R$. We will also use the shifted action 
$w \cdot \mu = w(\mu +\rho) -\rho$ for any $w \in W$, $\mu \in P$, where 
$\rho = \frac{1}{2} \sum_{\alpha \in R_+} \alpha$.  
Let $\Wl$ denote the affine Weyl group generated by the reflections 
$\{s_{\alpha, k} \}$ for any $\alpha \in R_+$, $k \in {\mathbb Z}$. The natural 
and shifted actions of $\Wl$ are defined by 
$s_{\alpha, k} ( \mu )= s_{\alpha} ( \mu) + kl\alpha$ and 
$s_{\alpha, k} \cdot \mu = s_{\alpha} \cdot \mu + kl\alpha$ for any $\mu \in P$.
The sets  $ \hat{X}= \{ \lambda \in P_+ : \langle \lambda , \alpha \rangle \leq  l \;\;{\rm for}\;{\rm all}\; \alpha \in R_+ \} $
and
$ \bar{X}= \{ \lambda \in P : 0\leq \langle \lambda + \rho, \alpha \rangle \leq  l \;\;{\rm for}\;{\rm all}\; \alpha \in R_+ \} $ 
are the fundamental domains for the non-shifted and shifted affine Weyl group actions respectively. Since the group of translations by $lP$ is also normalized by $W$, one can consider the semidirect product $\WlP \simeq W \ltimes lP$, which also naturally acts on $P$. 
We define the shifted action of $\WlP$ on $P$ by $w \cdot \lambda = w(\lambda +\rho) -\rho$ for 
any $w \in \WlP$, $\lambda \in P$. 

 Let $\Pi \in R_+$ denote the simple roots, and set $r = |\Pi|$, the rank of $\Lg$.   The affine Weyl group $\Wl$ can be viewed as a Coxeter system $(\Wl, S)$, where 
$S = \{ s_0, s_1, \ldots, s_r \}$ is the set of reflections with respect to the walls of $\bar{X}$. Here $\{s_1, \ldots s_r \}$ are the reflections in $W$ corresponding to the simple roots, 
and $s_0$ is the reflection about the upper wall of $\hat{X}$. The length $l(w)$ of an element in $\Wl$ is defined as the smallest number $n$, such that $w$ is a product of $n$ elements of  $S$.   

Let $\Ug$ be the Drinfeld-Jimbo quantum enveloping algebra associated to a semisimple simply-laced Lie algebra $\Lg$ over ${\mathbb Q}(v)$, where $v$ is a formal variable. It is 
generated by the elements $\{ E_i, F_i, K_i ^{\pm 1}\}_{i=1}^r$,  with the standard set of 
relations (see e.g. \cite{Lus1}). Denote by $U_v(\Lg)^{\mathbb Z}$ 
the ${\mathbb Z}[v,v^{-1}]$ subalgebra of $\Ug$ generated by the divided 
powers $E_i^{(k)} = \frac{E_i^k}{[k]!}$, $F_i^{(k)} = \frac{F_i^k}{[k]!}$  and $K_i, K_i^{-1}$ for $i =1 \ldots r$,
$k \in {\mathbb N}, k \geq 1$, where
 $[k]! = \prod_{s=1}^k \frac{v^s -v^{-s}}{v -v^{-1}}.$
Then the restricted quantum group $\Uq$ is defined by 
$\Uq = U_v(\Lg)^{\mathbb Z} \otimes_{{\mathbb Z}[v, v^{-1}]} {\mathbb Q}(q)$, 
where $ v$ maps to $q \in {\mathbb C}$, which is set to be a primitive $l$-th root of unity. 
By \cite{Lus1}, we have for $1 \leq i \leq r$: $E_i^l =0$, $F_i^l =0$, 
$K_i^{2l} =1$, and $K_i^l$ is central in $\Uq$. 

Let  $\Uf$  be the finite dimensional subquotient of $\Uq$ \cite{Lus2}, generated by 
$\{ E_i, F_i, K_i^{\pm 1} \}_{i=1}^r$ over ${\mathbb Q}(q)$, factorized 
over the two-sided ideal $\langle \{K_i^l -1\}_{i=1}^r \rangle $. Then  $
\Uf$  is a Hopf 
algebra of dimension $l^{{\rm dim}\Lg}$ over ${\mathbb Q}(q)$. 
We will use the same notations for the $\mathbb C$-algebras 
$\Uq$ and $\Uf$, where the 
field ${\mathbb Q}(q)$ is extended to $\mathbb C$.  

Let $\Cf$ denote the category of finite dimensional modules of type $\bf 1$ over $\Uq$, i.e. such that the central elements $\{K_i^l \}_{i =1\ldots r}$ act by $1$ in any object of $\Cf$.
Let $\bCf$ be the category of finite dimensional modules over $\Uf$. 
For a module $M$ in $\Cf$ (resp. $\bCf$) we will denote by $[M]$ its image in the Grothendieck ring $[\Cf]$ (resp. $[\bCf]$). We denote by $\mR$ (resp. $\R$) the complexification of the ring $[\Cf]$ (resp. $[\bCf]$). We will use the symbol $chM$ for the formal character 
$chM = \sum_{\eta \in P} ({\rm dim} M_\eta) e^\eta \in {\mathbb Z}[P]^W$, 
where we write 
$e^\eta$ for the basis element in ${\mathbb C}[P]$ corresponding to 
$\eta \in P$. Then ${\mathbb C}[P]$ can be considered 
as a space of functionals on the subalgebra $(\Uq)^0 \subset \Uq$, which is
generated by $K_i^{\pm 1}$ and the products 
$\prod_{s=1}^l \frac{K_i q^{1-s} - K_i^{-1} q^{s-1}}{q^s - q^{-s}}$ for 
all $i = 1, \ldots, r$; in particular 
$e^\eta(K_\beta) =q^{<\eta, \beta>}$ for any $\eta \in P, \beta \in Q$. 
Throughout the paper we will denote by $f(\lambda)$ the 
normalized $W$-symmetric function 
$\frac{1}{|W_\lambda|}\sum_{w \in W} e^{w(\lambda)} \in {\mathbb C}[P]^W$, 
where $W_\lambda$ is the stabilizer subgoup of $\lambda$. The Weyl character 
will be denoted by $\chi(\lambda)$ with $\chi(w \cdot \lambda) \equiv 
\varepsilon(w) \chi(\lambda)$, $\lambda \in P$. 
We will use the 
 symbol $[f(\lambda)]$ for the element of $\mR$ or $\R$, 
mapped to $f(\lambda)$ under the 
isomorphism $[M] \to ch M$. 

\subsection*{Acknowledgements} 
I am very grateful to my advisor Igor Frenkel for his guidance 
and valuable suggestions, and to Viktor Ostrik for many hepful 
discussions. 
I also would like to thank Aleksander Kirillov-Jr. and 
 the referee of the first version of this 
paper for the important corrections and suggestions. 
I am grateful to the Mathematics Department of Yale University 
where this work was done, and the Clay Mathematics Institute for the 
financial support.

\section{The fusion category} \label{sec:fus}

Inequivalent simple finite dimensional modules in $\Cf$ are in one-to-one correspondence with the dominant weights $P_+$ \cite{Lus1}. 
The Weyl modules, also parametrized by $P_+$, can be defined as the 
quantum deformations of simple finite dimensional modules for the 
classical universal enveloping algebra $U(\Lg)$ \cite{Lus1}. 
They play the role of the universal highest weight modules in 
$\Cf$. Both the Weyl modules and their duals have their characters given by the Weyl formula, though in general they are not simple over $\Uq$. It follows that the Grothendieck ring of the non-semisimple category $\Cf$ is isomorphic to  
${\mathbb Z}[P]^W$, the isomorphism given by $[M] \to chM$. The category $\Cf$ can be decomposed into a direct sum of subcategories according to the {\it linkage principle} \cite{APW1}: 
If simple modules $L(\mu_1)$ and $L(\mu_2)$ both are composition factors of an indecomposable module $M \in \Cf$, then $\mu_1 \in \Wl \cdot \mu_2$. 
A {\it tilting module} is defined as an element $T$ of $\Cf$ such that both $T$ and its dual $T^*$ have filtrations with factors isomorphic to  Weyl 
modules  \cite{And}, \cite{APW1}. The indecomposable tilting modules are 
also parametrized by $P_+$. One of the most important properties of the 
non-abelian category of tilting modules $\mathcal T$ is that it is closed 
with respect to tensor products \cite{And}, \cite{Par}. In the present work 
we will need the description of a semisimple quotient of $\mathcal T$, 
the {\it fusion category}, which we briefly outline below, 
following \cite{AP}. 

\subsection*{The fusion category $\mathcal Fus$ over $\Uq$.}

Let $\mR \simeq {\mathbb C}[P]^W$ denote the algebra of characters of the category $\Cf$, which is the complexification of the Grothendieck ring $[\Cf]$. 
Then $\mR$ has a natural basis of simple characters $\{[L(\lambda)]\}_{\lambda \in P_+}$, as well as a basis of Weyl characters $\{[W(\lambda)]\}_{\lambda \in P_+}$, and of indecomposable tilting characters $\{[T(\lambda)]\}_{\lambda \in P_+}$. 

Define the {\it first dominant alcove} 
$$ X = \{ \lambda \in P_+ : \langle \lambda + \rho, \alpha \rangle < l 
\;\;{\rm for}\;{\rm all}\; \alpha \in R_+ \}. $$

 \begin{Thm} \cite{AP} \label{thm:J} 
Decompose the algebra of characters of the category $\Cf$ into a sum 
$$ \mR = Vr \oplus J, $$
where 
$$ Vr  = {\rm span}_{\mathbb C}\{[T(\lambda)], \lambda \in X \},$$
$$ J = {\rm span}_{\mathbb C}\{[T(\lambda)], \lambda \in P_+\setminus X \}.$$
Then 
$$J = {\rm span}_{\mathbb C}\{ [W(\lambda)] +[W(s\cdot \lambda)] \}_{\lambda, s\cdot \lambda 
\in P_+ {\rm are}\; {\rm in}\; {\rm adjacent}\; {\rm alcoves}, s - {\rm a}\; {\rm refl.} \; {\rm in}\, \Wl}. $$
\end{Thm}

We will need the following consequence of this result. 
Define the elements $[W(\lambda)]$ for non-dominant weights by setting  
$[W(w\cdot \lambda)] = \varepsilon(w) [W(\lambda)]$ for any $\lambda \in P_+$, $w \in W$.

\begin{Cor} \label{cor:J}
In the notation of \thm{J},
$$J = {\rm span}_{\mathbb C}\{ [W(\lambda)] +[W(s\cdot \lambda)] \}_{\lambda \in P, s - {\rm a}\; {\rm reflection} \; {\rm in}\, \Wl}. $$
\end{Cor}

\begin{proof}
  Suppose that $s_H \in \Wl$ is a reflection about a hyperplane $H \subset P$,
$\lambda ' = s_H \cdot \lambda$ and $\lambda', \lambda$ are not in the adjacent alcoves. Then $s_H$ can be written as a composition of reflections
$s_H = s_1 \ldots s_k$ such that $s_1 \cdot (\lambda_0=\lambda) = 
\lambda_1, \ldots ,$  
$s_k \cdot \lambda_{k-1} = \lambda_k = \lambda'$ and each two weights 
$\lambda_i, \lambda_{i+1}$ lie in adjacent alcoves. The number of reflections 
$k$ is odd, and therefore any 
$[W(\lambda)]+[W(s_H \cdot \lambda)] = \sum_{i=0}^{k-1} (-1)^{i} 
([W(\lambda_i)]+[W(s_{i+1} \cdot \lambda_i)]) $ 
belongs to $J$. By the definition of $[W(\lambda)] \in \mR$ for $\lambda \notin P_+$, we have  $[W(\lambda)]+[W(s_i \cdot \lambda)]=0$ if $s_i$  is a reflection about one of the hyperplanes through $-\rho$, then the condition that $\lambda, s\cdot \lambda \in P_+$ can be dropped, 
and we get the statement of the Corollary.
\end{proof}

The subspace $J$ is a multiplicative ideal in $\mR$ (\cite{AP}, 3.19) and therefore the quotient   
$Vr \simeq \mR/J$ is a commutative associative algebra with a unit (the Verlinde algebra). Multiplication rules in $Vr$ can be written in the basis of tilting characters:
$$ [T(\lambda)] [T(\mu)] = \sum_{\nu \in X} a_{\lambda, \mu}^\nu [T(\nu)] 
\,({\rm mod}J),$$ 
where $a_{\lambda, \mu}^\nu$, $\lambda, \mu \in X$, $\nu \in P_+$ 
are the coefficients of decomposition of the tensor product of two tilting 
modules into indecomposables. For $\lambda \in X$, each tilting module is a Weyl module and its 
character is given by a unique element $[W(\lambda)] \in \mR$.

Let ${\mathcal T}'$ be the subcategory of $\mathcal T$ whose objects are tilting modules with no indecomposable component of highest weight in $X$.
Then by the above, 
the subcategory ${\mathcal T}'$ of modules with characters in $J$ is a tensor ideal in 
$\mathcal T$. This gives rise to a reduced tensor product structure of the   category 
${\mathcal Fus} \simeq {\mathcal T}/{\mathcal T}'$ with the fusion ring defined above. 

The ideal ${\mathcal T}' \subset {\mathcal T}$ admits an equivalent 
characterization in terms of the quantum traces. 
For $\mu = \sum_{i=1}^r 
n_i \alpha_i \in Q$, $\alpha_i \in \Pi$, define 
$K_\mu \equiv  
\prod_{i=1}^r K_i^{n_i}$, where $n_i \in {\mathbb Z}$. For any $M 
\in \Cf$ and any $\Uq$-endomorphism $\phi$ of $M$, define the quantum 
trace $tr_q(\phi)$ by setting $tr_q(\phi) = tr(K_{2\rho}\phi)$, 
where $K_{2\rho} = \prod_{\alpha \in R_+} K_{\alpha}$ and 
$tr$ is the usual trace on $M$. The quantum dimension of a module is defined to be 
the quantum trace of the identity morphism. 

\begin{Prop} \cite{And} \label{prop:qtrace} 
Let $M \in {\mathcal T}'$. Then for any $\Uq$-endomorphism $\phi$ of $M$, 
$tr_q(\phi)=0$.
In particular, $dim_q(M) =0$.
\end{Prop}   

Now we will consider the corresponding categorical and algebraic constructions over the small quantum group $\Uf$. 

\subsection*{Restriction of $\Cf$ to $\Uf$.}

The category $\bCf$ of finite dimensional modules over $\Uf$ inherits many 
of its properties from the category $\Cf$ \cite{Lus2}, \cite{APW2}.
Simple modules in $\bCf$ are parametrized by their highest weights, and inequivalent simple modules are the restrictions of those over $\Uq$ with the restricted highest weights $P_l 
\equiv \{ \lambda \in P_+ : \langle \lambda, \alpha_i \rangle \leq l-1 \; {\rm for}\; 
{\rm all}\; \alpha_i \in \Pi\}$. 
Define the character of a simple module $L(\lambda)$, $\lambda \in P_l$ over 
$\Uf$ to be the corresponding element $[L(\eta)]$ of the Grothendieck 
ring $[\bCf]$ of the category 
$\bCf$. Then for any finite dimensional $\Uf$-module $M$ with composition 
factors $L(\eta)$ of multiplicities $n_\eta$, define the character of $M$ 
by the formula 
$$ [M] = \sum_{\eta \in P_l} n_\eta [L(\eta)] \in [\bCf].$$
The characters of simple modules $\{[L(\eta)]\}_{\eta \in P_l}$ form a basis 
in $\R$, the complexification of the Grothendieck ring $[\bCf]$. 
The multiplicative structure of the algebra $\R$ 
 can be obtained from that of $\mR \simeq {\mathbb C}[P]^W$ using the 
tensor product decomposition for  
simple modules $\{L(\mu) \}_{\mu \in P_+}$ over $\Uq$:
\begin{Thm} \cite{Lus1} \label{thm:fs} 
Let $\mu = \mu_0 + l\mu_1 \in P_+$ be a dominant weight such that 
$0 \leq \langle \mu_0, \alpha_i \rangle \leq l-1$ for all $\alpha_i \in \Pi$. Then  
$$ L(\mu_0 + l\mu_1) \simeq L(\mu_0) \otimes L(l\mu_1)$$ 
as $\Uq$-modules. The restriction of $L(\mu_0)$ to $\Uf$ is a simple 
$\Uf$-module. 
\end{Thm} 

The subquotient  $\Uf$  of $\Uq$ acts trivially on the simple module $L(l\mu_1)$, 
whose weights are all in $lP$. By restriction, simple $\Uq$-modules of type $\bf 1$ are  
modules over $\Uf$. Therefore, the algebra of characters 
$\R$ of  $\Uf$-modules  can be obtained 
from that of  $\Uq$  by setting the characters of simple modules with $l$-multiple highest weights equal to their dimensions, and we get
\begin{Cor} \label{cor:R}
The following diagram commutes: 
$$ \begin{CD} \mR @>>> \R  \\
@V{\wr}VV             @V{\wr}VV   \\
{\mathbb C}[P]^W  @>>> {\mathbb C}[P]^W \otimes_{{\mathbb C}[lP]^W} {\mathbb C},\end{CD} $$
where $ {\mathbb C}[lP]^W \longrightarrow {\mathbb C}$ maps each element to its value at $1$. 
\end{Cor}

The isomorphisms $\mR \simeq {\mathbb C}[P]^W$ and $\R \simeq 
({\mathbb C}[P]^W \otimes_{{\mathbb C}[lP]^W} {\mathbb C})$ of \cor{R} will 
be denoted by $\psi$ and $\bar{\psi}$,  respectively. For 
an element $g \in \mR$, we will write $\bar{g}$ for its image in $\R$ 
under the 
homomorphism $\mR \to \R$ induced by the restriction of $\Uq$-modules to 
$\Uf$. 

We are interested in the multiplicative structures associated to the category 
$\bCf$. Since it is not semisimple (which is inherited from non-semisimplicity of $\Cf$), the algebra of characters described in \cor{R}
gives only partial information, namely the multiplicities of simple 
modules in the filtration of a tensor product. As in the case of ${\mathcal T} \subset \Cf$, 
it would be convenient to define a subcategory in $\bCf$ designed to behave nicely with respect to tensor products. This suggests to consider the restriction of 
tilting modules over $\Uf$. 

\begin{Prop} 
The characters $\{ [T(\lambda)]\}_{\lambda \in P_l}$ form 
a basis in $\R$. 
\end{Prop} 

\begin{proof} Define a partial ordering in $P_l$ by $\lambda 
{\geq}^l \mu$ if 
there exists a weight $\mu_1 \in P_+$ such that $\lambda \geq \mu + l\mu_1$ in 
the usual ordering in $P$. Then the basis of tilting characters is 
connected to $\{[L(\eta)]\}_{\eta \in P_l} \in \R$ 
by a triangular transformation with respect to the ordering 
${\geq}^l$. 
\end{proof}

It should be mentioned that the restrictions of tilting modules are not 
tilting over $\Uf$ in the sense of the abstract definition. Still, since 
$\Uf$ inherits the coalgebra structure from $\Uq$, the 
set of restrictions of tilting modules to $\Uf$ is closed under 
the tensor product. 

\subsection*{The fusion category considered over $\Uf$.}

The fusion category was constructed over $\Uq$  as a quotient of the category of tilting modules. The same structure can be considered over the finite dimensional algebra $\Uf$. Let 
$\bar{\mathcal T}$ be the subcategory of $\bCf$ whose objects are the restrictions of $\Uq$-tilting modules to $\Uf$. Then it is closed under the tensor products.  The ideal 
${\mathcal T}' \subset {\mathcal T}$ restricts to the set of the restrictions 
of  tilting modules with indecomposable components (over $\Uq$) of highest 
weights in $P_+ \setminus X$. The subcategory $\bar{\mathcal T}'$ of $\bar{\mathcal T}$ whose objects are the restrictions of tilting modules in 
${\mathcal T}'$, is closed under tensor products. 

Let $\bVr$ denote the quotient ${\R}/{\bJ}$, where $\bJ$ is the image of 
the ideal $J$ under the homomorphism $\mR \to \R$, sending the character 
of a $\Uq$-module to the character of its restriction over $\Uf$. 
Then the algebra structure in $\bVr$ corresponds to the truncated tensor 
product of the restrictions of tilting modules with no indecomposable 
components of highest weight in $P_+ \setminus X$.

To describe the ideal $\bJ$ of $\R$, we have to introduce the action of 
the affine Weyl group in $P_l$.  

Note that the natural (non-shifted) action of the affine Weyl group $\Wl$ on $P$ preserves the $lP$-lattice in $P$:
$$ s_{\alpha, k} (l\mu) = l\mu - l \langle \alpha,  \mu \rangle \alpha + k l \alpha \in lP.$$ 
Similarly, the shifted action of $\Wl \cdot$ on $P$ preserves the shifted lattice $lP - \rho$. 
Therefore, the following actions of $\Wl$ on $P_l$ can be defined: 
\begin{enumerate}
\item{the natural (non-shifted) action
$$ w \circ \lambda = w(\lambda) \,{\rm mod}(lP)$$ 
for any $\lambda \in (P/lP)_{+}$, $w \in \Wl$, and}
\item{the shifted action 
$$ w \bullet \lambda = (w \cdot \lambda) \,{\rm mod}(lP)$$
for any $\lambda \in P_l$, $w \in \Wl$}
\end{enumerate} 
By definition, the results of both actions belong again to $P_l$.

Since $lQ \subset lP$, it is enough to consider the $\circ$- and 
$\bullet$-actions of elements $w \in W$ of the finite Weyl group.

Recall that the affine Weyl group $\Wl$ is isomorphic to $W\ltimes lQ$, where $lQ$ generate the translations by $l$-multiple roots, and consider the affine group $\WlP = W \ltimes lP$, where the translations are generated by the $lP$ lattice. Let $\Omega$ be the subgroup of $\WlP$ stabilizing the alcove $\bar{X}$ (resp. $\hat{X}$) with respect to the shifted (resp. non-shifted) action. Suppose that $lP \cap Q = lQ$. This 
condition is always satisfied if the number $l$ is such that ${\rm GCD}(l, {\rm det}a_{ij})=1$ for the Cartan matrix $a_{ij}$ of $\Lg$ (\cite{Lyu1}, \cite{La}). 
To simplify the subsequent considerations, we will assume that this condition 
holds, and also that
$l \geq h$ for the root system $R$ of $\Lg$.
Then by \cite{Hum3} $\Wl \cap \Omega =1$, so in fact $\WlP$ is the semidirect product of $\Wl$ and $\Omega$, and $\Omega \simeq \WlP/\Wl$. The order of $\Omega$ is equal to $|\WlP/\Wl|=|P/Q|$, the index of connection. 

Denote by $\hat{\mathcal X}$ the set of representatives for the orbits 
of the 
non-shifted action of $\Omega$ in the alcove $\hat{X}$. Similarly, let 
$\tmX = \bar{X}/\Omega$ be the set of representatives for 
the orbits of $\Omega$ with respect to the shifted action. Then $\hat{\mathcal X}$ and $\tmX$ are the fundamental domains for  the $W \circ$ and $W \bullet$-actions respectively, and $|\hat{\mX}| = |\tmX|$. 
Let $\mathcal X \subset \tmX$ denote the subset of regular weights in $\tmX$, such that their stabilizer in $W$ with respect to the $\bullet$-action is trivial. 

The $W \bullet$-action is designed to formulate the linkage principle for $\Uf$:
\begin{Prop} \label{prop:link}
 If two simple modules $L(\mu)$ and $L(\lambda)$ both are composition factors in a filtration of an indecomposable $\Uf$-module, then $\mu \in W \bullet \lambda$.
\end{Prop}
\begin{proof}

Consider the graded version $\bUf$ of $\Uf$ with $(\bUf)^0 \simeq (\Uq)^0$,
as defined 
in e.g. \cite{AJS}, \S 1.3, in which the toric part coincides with that 
of $\Uq$. 
Then two simple modules 
$\stackrel{\bullet}{L}(\lambda)$ and $\stackrel{\bullet}{L}(\mu)$ both 
are composition factors in an indecomposable $\bUf$-module
if and only if $\lambda \in \Wl \cdot \mu$ \cite{AJS}, \S 6.7. 
In the non-graded setting this means that any two simple composition factors 
of an indecomposable module have their    
 highest weights in a 
$\Wl \cdot$-orbit restricted modulo $lP$, which coincides with the 
definition of an orbit of $W \bullet$-action. In particular, this means 
that the set $\tmX$ enumerates the blocks of the category $\bCf$.
\end{proof}

\begin{Thm} \label{thm:bJ}
Let $\bJ \subset \R$ be the image of $J \subset \mR$ under the 
homomorphism $\mR \to \R$. Then 
$$\bJ = span_{\mathbb C}\{[W(\lambda)] + 
[W(s\bullet \lambda)] \}_{\lambda \in P_l, s - {\rm a}\; 
{\rm reflection} \; {\rm in}\, W}. $$
\end{Thm} 

\begin{proof}

Use the notations introduced after \cor{R}. 

First we will show that for any $\lambda \in P_l$ and a reflection $s \in W$,
the element $[W(\lambda)] + [W(s\bullet \lambda)]\equiv \bar{g}_\lambda$ is 
in $\bJ$. 

Decomposing the Weyl characters in terms of simple composition factors over 
$\Uf$, we can write $\bar{g}_\lambda = \sum_{w \in W}a_w [L(w \bullet \lambda)]$ for 
some coefficients $a_w \in {\mathbb C}$. 
For each $w \in W$ the preimage of $[L(w \bullet \lambda)]$ under the 
homomorphism $\mR \to \R$  contains 
a scalar multiple of an element of the form $[L(z \cdot \lambda)]$ 
for some $z \in \Wl$.   
Then $g_\lambda$ can be written as $g_\lambda = \sum_{z \in \Wl}b_z 
[L(z\cdot \lambda)]$, or in the basis of tilting characters, 
$g_\lambda = \sum_{z \in \Wl}c_z [T(z\cdot \lambda)]$ 
for some coefficients $b_z, c_z \in {\mathbb C}$. 
We have 
$$\begin{array}{c}
\bar{\psi}(\bar{g}_\lambda)(K_{2\rho})= (\chi(\lambda)+ \chi(s\bullet \lambda))(K_{2\rho}) = \\ 
=  \frac{1}{\Delta(K_{2\rho})} 
\sum_{w \in W}(\varepsilon(w)+\varepsilon(ws^{-1})) 
q^{\langle w(\lambda+\rho),2\rho \rangle} = 0,
\end{array}$$ 
where $\Delta$ denotes the Weyl denominator. 
For any $g \in \mR$ and $\mu \in P_l$, we have $\psi(g)(K_{\mu})= \bar{\psi}(\bar{g})(K_\mu)$, and therefore 
$$ \psi(g_\lambda)(K_{2\rho}) = \sum_{z \in \Wl}c_z chT(z \cdot \lambda)(K_{2\rho}) =0. $$ 
If $\lambda$ is singular with respect to the $\Wl \cdot$ action, then 
the identity always holds. If $\lambda$ is regular, then it requires $c_1 =0$. 
In both cases $g_\lambda \in J$, and therefore 
$\bar{g}_\lambda \in \bJ$. 

Now we will show that $\bJ$ cannot be bigger than required in the Theorem. 
Choose representatives for regular $W\bullet$-orbits in $P_l$ so that 
$\mX \subset X$, and 
consider the characters $[T(\lambda)]=[L(\lambda)]$ for $\lambda\in \mX$. 
Suppose that for some coefficients $d_\lambda \in {\mathbb C}$ there exists 
a linear combination $ \bar{t} \equiv \sum_{\lambda \in \mX} 
d_\lambda [T(\lambda)] \in \R$ of characters of restrictions of tilting 
modules, which belongs to $\bJ$. Then its pre-image under the homomorphism 
$\mR \to \R$ lies in $J$. Since  
$[T(\lambda)]_{\lambda \in \mX}$ correspond to pairwise distinct 
$W\bullet$-orbits 
in $P_l$, each of their pre-images $t_\lambda$ lies in $J$. Then 
$d_\lambda \psi(t_\lambda) (K_{2\rho})=0$ 
for each $\lambda \in \mX$, and hence $d_\lambda ch T_\lambda (K_{2\rho})=0$,
which by construction of $\{[T(\lambda)]\}_{\lambda \in \mX}$ requires each 
$d_\lambda =0$. Therefore, no linear combination of characters of the 
restrictions of $\{T(\lambda)\}_{\lambda \in \mX}$ lies in $\bJ$, and 
${\text dim} \R/\bJ \geq |\mX|$. On the other hand, 
$${\text dim} \R/\langle [W(\lambda)] + 
[W(s\bullet \lambda)] \rangle_{\lambda \in P_l, s - {\rm a}\; 
{\rm reflection} \; {\rm in}\, W} = |\mX |, $$
and the theorem follows. 
\end{proof}

 The dimension of $\bVr$ coincides with the number of $W\bullet$-regular 
weights $\mX \subset X$ in $\tmX$, and the elements  
$\{[W(\lambda)] =[L(\lambda)]\}_{\lambda \in {\mX}}$ form a basis in $\bVr$. 

With the assumption on the number $l$ introduced above, the order $|X|$ is divisible by $|P/Q|$. Then we get

\begin{Cor}

$${\rm dim}{\bVr} = {\rm dim}{Vr}/|P/Q|.$$
\end{Cor} 

\begin{Cor} \label{cor:Vr-}
$\bVr$ is isomorphic to the algebra 
$$ \bar{(Vr)}^- \simeq {\mathbb C}[P]^W/\langle (\chi(\lambda) + \chi(s \cdot \lambda)\rangle_{\lambda \in P, s - {\rm a}\,{\rm reflection}\, {\rm in}\, \WlP}.
$$
\end{Cor}
\begin{proof}
 Clearly, $\tmX$ is the fundamental domain for the 
shifted action of $\WlP$ on $P$, and $\mX \subset \tmX$ is the set of the regular weights with respect to this action.  The ideal $\bJ$ of $\R$ is the image 
of $J \subset \mR$ under the homomorphism $\mR \to \R$.
The ideal of $\mR$ spanned by $\{[W(\lambda)] + [W(s\cdot \lambda)] \}_{\lambda \in P, s \in \WlP, s^2=1}$ contains $J$, and therefore its image in $\R$ contains $\bJ$. Compare the dimensions: the elements $\{\chi(\lambda) \}_{\lambda \in \mX}$ form a basis in $ \bar{(Vr)}^-$. Therefore, 
${\rm dim} ( \bar{(Vr)}^-) = 
|\mX| = {\rm dim}(\R/\bJ)$, 
 which means that the algebras are isomorphic. 
\end{proof}

\section{Projective modules} \label{sec:Pr}

\subsection*{The category of projectives over $\Uq$.}

 Let $\mathcal Pr$ denote 
the category of projective modules over $\Uq$. Below we list their most important 
properties.

\begin{Prop} \cite{APW1} \label{prop:inj-proj}
\begin{enumerate} 
\item{The Steinberg module $St = L((l-1)\rho)$ is both projective and injective object in ${\mathcal C}_f$}
\item{All injectives in ${\mathcal C}_f$ are projectives and conversely.}
\item{All projective objects in ${\mathcal C}_f$ are tilting modules. Every 
indecomposable projective module is isomorphic to some $T(\lambda)$, 
$\lambda \in (l-1)\rho + P_+$.}
\item{Any indecomposable projective module in ${\mathcal C}_f$ is a direct 
summand of $St \otimes E$ for some $E \in {\mathcal C}_f$.} 
\end{enumerate} \end{Prop}

Therefore, the set ${\mathcal Pr} \subset \Cf$ is in fact a subcategory in 
$\mathcal T$. It is also a tensor ideal in $\mathcal T$, since a direct 
summand of a projective module is projective, and a tensor product of 
a projective module with any module in $\Cf$ is again projective. 
Since $St$ is a Weyl module, its character $ch St = \chi((l-1)\rho) 
\in {\mathbb C}[P]^W$ 
is given by the Weyl character formula. 

\begin{Prop} \label{prop:proj_ch}
The linear span of the characters of projective (injective) modules in 
${\mathbb C}[P]^W$ coincides with $ch St \cdot {\mathbb C}[P]^W$. 
\end{Prop}
\begin{proof}
 Any $p \in {\mathbb C}[P]^W$ can be written as a linear combination of characters of simple modules in ${\mathcal C}_f$. Since $St$ is projective, any element of the form $ch St \cdot p \in  {\mathbb C}[P]^W$ is a linear combination of characters of projective modules.
 Then we have  
$$ \{ ch St \cdot f(\lambda) \}_{\lambda \in P_+} \subset 
span_{\mathbb C}\{ch P(\mu) \}_{\mu \in (l-1)\rho +P_+}, $$ 
where $f(\lambda)$ is the $W$-symmetric function (\S 1), and $chP(\mu)$ denotes the image in 
${\mathbb C}[P]^W$ of the character of the projective module 
$P(\mu) = T(\mu)$ with highest weight $\mu$. Note that in the decomposition of $ch St \cdot f(\lambda)$ into a linear combination of projective characters, only the highest weights $\mu : (l-1)\rho \leq \mu \leq (l-1)\rho +\lambda$ can occur. 
Therefore, the two sets are related by a lower triangular transformation of bases, which provides a bijection between the span of indecomposable projective characters and  $ch St \cdot {\mathbb C}[P]^W$. 
\end{proof}

The next proposition provides the factorization of projective (injective)
modules over $\Uq$ similar to the statement of \thm{fs}.

\begin{Prop} \cite{And} 
\label{prop:fp} 
For $\lambda = \lambda_0 + l \lambda_1 \in P_+$, where 
$0 \leq \langle \lambda_0, \alpha_i \rangle \leq l-1, \; i=1, \ldots, r$,
let $\bar{\lambda}$ denote the weight $\bar{\lambda} =
2(l-1)\rho + w_0 (\lambda_0) + l \lambda_1 \in P_+$. Then 
\begin{enumerate} 
\item{The indecomposable tilting module $T(\bar{\lambda})$ is both projective 
and injective in $\Cf$;}
\item{$T(\bar{\lambda}) \simeq T(\bar{\lambda_0}) \otimes L(l\lambda_1).$}
\end{enumerate} 
\end{Prop} 

Now we are can describe the ideal $Pr \in \mR$.
For any $\mu \in P$, let $\hat{W}_\mu$ denote the subgroup 
of $\Wl$ generated by 
$\{s_{\alpha, \langle \mu, \alpha \rangle} \}_{\alpha \in \Pi}$. Then 
$\hat{W}_\mu$ stabilizes $l\mu -\rho$ with respect to the shifted action,
and is isomorphic to $W \subset \Wl$. 

\begin{Thm} \label{thm:Pr}
The characters of projective modules over 
$\Uq$ span the ideal of $\mR$ of the form      
$$ Pr = span_{\mathbb C}\{ \sum_{w \in \hat{W}_\mu}[W(w\cdot \lambda)]\}_{\lambda 
\in P, \mu \in P}. $$
\end{Thm} 

\begin{proof}

We will use 
the isomorphism $\psi : \mR \xrightarrow[\sim]{} {\mathbb C}[P]^W$ and prove the corresponding  
statement for $\psi(Pr) \subset {\mathbb C}[P]^W$.  
The ideal of projective characters coincides with the linear span of 
$\{ch St \cdot f(\lambda)\}_{\lambda \in P_+} = \{ \sum_{w \in \hat{W}_\rho} 
\chi(w\cdot((l-1)\rho + \lambda)\}_{\lambda \in P_+}$, thus all projective characters are in $\psi(Pr)$.
On the other hand, for any $\lambda, \mu \in P$, the element  
$$ \sum_{w \in \hat{W}\mu}\chi(w\cdot \lambda) $$
is either zero, if $\langle \mu, \beta \rangle = 0$ for some 
$\beta \in R$, or is equal to 
$$ \chi((l-1)\rho + l(\mu-\rho))\cdot f(\lambda - l\mu). $$ 
 If $(\mu -\rho) \in P_+$, 
$$ W((l-1)\rho + l (\mu - \rho)) \simeq St \otimes L(l(\mu -\rho))$$  
is a projective and simple module \cite{APW1}, and therefore 
$ \chi((l-1)\rho + l(\mu - \rho))\cdot f(\eta)$ lies in the linear span 
of projective characters in ${\mathbb C}[P]^W$ for any $\eta \in P$.
For nondominant weights, $\chi((l-1)\rho + l(\mu -\rho))$ is either zero or 
equals   
$\pm \chi((l-1)\rho + l\nu)$ for some $\nu \in P_+$, and the Proposition follows. 
  
\end{proof} 

Thus, both $J$ and $Pr$ admit characterizations in terms of the 
symmetries with respect to the shifted action of the affine 
Weyl group (cf. \cor{J}), and $Pr \subset J$. 

\subsection*{Projective modules over $\Uf$}

Next we are interested in the 
restrictions of projective modules to $\Uf$. 
Below we formulate the known facts on projective modules in $\bCf$. 

\begin{Thm} \cite{APW2} \label{thm:rp} 
\begin{enumerate}
\item{ $St$ is both projective and injective in $\bCf$.}
\item{All injective modules in $\bCf$ are projective, and vice versa.}
\item{Keep the notations of \prop{fp}.
For each $\lambda \in P_+$, the restriction of $T(\bar{\lambda})$ 
to $\Uf$ is injective and projective in $\bCf$.}
\item{For each $\lambda_0 \in P_l$, the restriction of 
$T(\bar{\lambda_0})$ 
to $\Uf$ is isomorphic to the injective hull of $L(\lambda_0)$ in $\bCf$. 
In particular, any injective $\Uf$-module is the restriction of an injective $\Uq$-module.}
\item{For each $\lambda \in P_+$, the restriction of $T(\bar{\lambda})$ 
to $\Uf$ is isomorphic to the direct sum of $dim(L(l\lambda_1))$ copies of 
$T(\bar{\lambda_0})$.}
\end{enumerate} \end{Thm}

Projective modules form a tensor ideal $\bar{\mathcal Pr}$ in the category of the restrictions of tilting modules to $\Uf$, and the linear span of their characters is a multiplicative ideal in $\R$. 
\begin{Prop} \label{prop:pr_f}
Linear span of the characters of projective (injective) modules over $\Uf$ coincides with 
$$ ch St \cdot {\mathbb C}[P]^W \otimes_{{\mathbb C}[lP]^W} {\mathbb C}.$$ 
\end{Prop} 
\begin{proof}
This follows from \prop{proj_ch} and the factorization of 
${\mathbb C}[P]^W$ to the restricted weights (\cor{R}).
\end{proof}

The next theorem determines the dimension of $\bPr \subset \R$:

\begin{Thm} \label{thm:bPr}
The characters of projective modules over 
$\Uf$ span the ideal of $\R$ of the form      
$$ \bPr = span_{\mathbb C}\{ [St]\cdot [f(\lambda)] \}_{\lambda \in \hat{\mX}}, $$ 
where $[f(\lambda)]$ denotes the element mapped to $f(\lambda)$ under the 
isomorphism of  \cor{R}, and $\hat{\mX}$ is the fundamental 
domain for the $W\circ$-action in $P_l$. 
\end{Thm}

\begin{proof}

\begin{Lem} \label{lem:propor}
For any $\lambda \in P_l$, the characters of any two projective $\Uf$-modules with no composition factors of highest weight outside the 
orbit $\{W \bullet \lambda \}$, are proportional 
in $\R$. 
\end{Lem}
\begin{proof}
  The radical of the algebra ${\mathbb C}[P]^W \otimes_{{\mathbb C}[lP]^W} {\mathbb C}$ is spanned by the elements which annihilate $K_\mu$ for any 
$\mu \in P_l$. It coincides with the linear span of the differences  
 $$ f(\lambda) - f(z \circ \lambda), $$
where $z \in W$, $\lambda \in P_l$ and we used the non-shifted $\circ$-action of $W$ on the restricted weights. 
 By \cite{La}, Lemma 5.3, the character of the Steinberg module 
$[St]$ annihilates the radical of $\R$, and  we have 
$$ chSt \cdot f(\lambda) - chSt \cdot  f(z \circ \lambda) =0 $$ 
for any $z \in W$, $\lambda \in P_l$.

By \prop{pr_f} any projective character with 
composition factors of highest weights in one orbit $W\bullet \lambda$ 
can be written as a linear combination 
$$ \sum_{z \in W} a_z chSt \cdot f(z \circ \lambda),$$ 
for some coefficients $a_z \in {\mathbb C}$, which by the above equals to 
$$(\sum_{z \in W} a_z) 
 chSt \cdot f(\lambda). $$
This shows that 
any two projective modules in the same block of the category $\bCf$ 
have equal characters up to a scalar multiple.  
\end{proof} 

By \thm{rp}, there exists at least one projective module in $\bCf$ with composition factors of highest weights in each $W\bullet$-orbit in $P_l$. Characters 
of modules in different blocks of category $\bCf$ are linearly independent 
in $\R$, and we get the statement of the Theorem. 
\end{proof} 

Note that by \lem{propor}, the ideal of characters of projective modules 
over $\Uf$ is $|\tmX|$-dimensional, where $\tmX$ is the fundamental domain 
for the $W\bullet$-action on $P_l$ ( we have $|\hat{\mX}| = |\tmX|$). Thus 
 $\bPr$ and $\bVr$ have bases parametrized respectively by the 
 weights in $\tmX$ and $\mX$, the closed and open alcove of the 
$W\bullet$-action in $P_l$.

\subsection*{Classification of tensor ideals in $\mathcal T$.}

 According to the main 
theorem in \cite{Ost}, the partially ordered (by inclusion) set of tensor ideals in 
$\mathcal T$ is isomorphic to the partially ordered set of two-sided cells in the 
affine Weyl group $\Wl$. The partition of a Coxeter group $G$ into a partially ordered set of two-sided cells was defined in \cite{KL} to study the ideals in the group algebra of $G$, consistent with a special (Kazhdan-Lusztig) basis in $G$. In particular, for the affine Weyl group $\Wl$, the unit element $\{1 \}$ is the highest, and the set $W_{(\nu)} \subset \Wl$ - the lowest cell with respect to this order \cite{Lus3}, \cite{Bed}, 
\cite{Shi}. 
The cell $W_{(\nu)}$ has the following description \cite{Shi}: 
$$  {W}_{(\nu)} = \{ w \in \Wl : w = x\cdot w_I \cdot y, \;{\rm for}\;
{\rm some}\; x,y \in \Wl,  \; l(w) = l(x) + l(w_I) + l(y) \},$$
where $w_I$ is the longest 
element of a standard parabolic subgroup $W^I \subset \Wl$, isomorphic to $W$. 

The correspondence in \cite{Ost} can be described as follows. Let $W^f$ be the set of minimal length representatives of the cosets $W\backslash \Wl$. Then by \cite{LX} 
for any two-sided cell $\underline{c} \in \Wl$, the intersection $\underline{c} \cap W^f$ is nonempty.  
Let $\underline{b} = \cup_i \underline{c_i} \in \Wl$ denote the union of all two-sided cells with the condition that $\underline{c_i} \leq \underline{c}$ for each $i$ 
with respect to the partial order defined in \cite{KL}. 
Then the full subcategory of $\mathcal T$ formed by the direct sums of 
tilting modules $T(\lambda)$ with $\lambda \in 
(\underline {b} \cap W^f)\cdot \bar{X}$ is  a tensor ideal in $\mathcal T$
corresponding to the cell  $\underline{c}$.  

\begin{Prop}
\begin{enumerate} 
\item{ The ideal ${\mathcal T}' \subset {\mathcal T}$ corresponds to the complement of 
the highest cell $(\Wl \setminus \{1 \}) \subset \Wl;$}
\item{  The ideal ${\mathcal Pr} \subset {\mathcal T}$ corresponds to the lowest cell 
$W_{(\nu)} \subset \Wl$.}
\end{enumerate} 
\end{Prop}
\begin{proof}
(1) was observed in \cite{Ost}: for $\lambda \in P_+$ we have $\lambda \in (W^f \setminus \{1\})\cdot \bar{X}$ if and only if $\lambda \in P_+\setminus X$,
which coincides with the set of highest weights of tilting modules in ${\mathcal T}'$.

For (2), we can use the description of $W_{(\nu)}$ in \cite{Bed}. For each pair of nearest parallel hyperplanes of affine reflections in P, define the strip 
$U$ to be the set of weights between these hyperplanes, and let $\mathcal U$ denote the set of all strips in $P$. Then by \cite{Bed}, 
$$ W_{(\nu)} \cdot \bar{X} = P\setminus (\cup_{U \in {\mathcal U} : 
X \subset U}U). $$
Therefore, $((W_{(\nu)} \cap W^f)\cdot \bar{X}) \cap P_+ = P_+ \cap ( P\setminus 
(\cup_{U \in {\mathcal U} : X \subset U}U))$, which coincides with the 
set 
$(l-1)\rho +P_+$. This is exactly the set of highest weights of the projective modules in $\mathcal T$.
\end{proof}

This means that $\mathcal Fus$ and $\mathcal Pr$ can be thought of as the 
"top" and the "bottom" of the tensor category $\mathcal T$. 
This relation also holds for their restrictions over $\Uf$:

\begin{Cor}
In the partially ordered (by inclusion) set of ideals in $\R$ corresponding 
to the restrictions of tensor 
ideals in $\mathcal T$, $\bJ$ and $\bPr$ are respectively the highest and 
the lowest nontrivial elements.
\end{Cor} 

Unlike in case of the big quantum group $\Uq$, here the ideal $\bPr$ is 
finite dimensional. The above argument then suggests that the ideal 
$\bPr$ should be viewed in the same framework with the  
 algebra $\bVr$, and, just as  
the Verlinde algebra, it should play a special role in the 
representation theory of $\Uf$. Further evidence of this similarity is 
presented in the next section.

\section{$\bPr$ as an analog of the Verlinde algebra.} \label{sec:VrPr}

In this section we will consider the properties of $\bPr$ in parallel to those of 
$Vr$ and its restriction to the small quantum group $\bVr$.
 
Write $\R = \oplus_{\mu \in \tmX} \R_\mu$ for the decomposition of 
the algebra $\R$ into a direct sum  of two-sided multiplicative 
ideals corresponding to the block decomposition of the category 
$\bCf$ (\prop{link}). The structure of 
the block $\R_\mu$ was described in \cite{La} following the result of \cite{BG}:
$$ \R_\mu \simeq {\mathbb C}[P]^{W_\mu} 
\otimes_{{\mathbb C}[P]^W} {\mathbb C}, $$
where $W_\mu$ is the subgroup of $W$ stabilizing the weight $\mu$ with respect to the $\bullet$- action, and ${\mathbb C}[P]^W \longrightarrow {\mathbb C}$ maps each $e^{l\eta}$ to $1$.

The ideal $\bPr$ has a clear interpretations in terms of the block 
structure of $\R$. The following statements can be derived from \thm{bPr} and 
\lem{propor}.

\begin{Cor} \label{cor:socR}
$$ \bPr \simeq {\rm Ann}({\rm Rad}\;\R).$$ 
\end{Cor}
\begin{proof}
The character of the Steinberg module annihilates the radical of $\R$ and is therefore an element of 
${\rm Ann}( {\rm Rad}\;\R)$. Since $[St]$ generates the algebra $\bPr$ over $\R$, any element of $\bPr$ belongs to ${\rm Ann}( {\rm Rad}\;\R)$. Now compare the dimensions: 
$$ |\tmX| = {\rm dim}(\bPr) \leq {\rm dim}{\rm Ann}({\rm Rad}\;\R) = |\tmX|,$$
where the first equality holds by \thm{bPr}, and the last follows from the block structure of $\R$.
Therefore, $\bPr \simeq {\rm Ann}({\rm Rad}\;\R)$.
\end{proof}

\begin{Cor} \label{cor:bPr}
The ideal $\bPr$ is isomorphic to $ [St] (\R/Rad(\R))$, where 
 the algebra $\R/{\rm Rad}(\R)$ is isomorphic to $(Vr)^+$ as defined in \S 1:
$$ (Vr)^+ \simeq {\mathbb C}[P]^W/\langle (f(\lambda) - f(s (\lambda))\rangle_{\lambda \in P, 
s - {\rm a}\,{\rm reflection}\, {\rm in}\, \WlP}$$ 
\end{Cor}
\begin{proof} Follows from the definition of $W\circ $ action on $P_l$.
Note that here in the definition of the ideal, a reflection $s \in \WlP$ can be replaced by any element 
$w \in \WlP$.  
\end{proof}

Choose the set of representatives $\hat{\mX}$ for the orbits of $W\circ$ in $P_l$ so that 
$\hat{\mX} \subset X$, by taking the smallest weight in each orbit. 
\begin{Prop} \label{prop:P(mu)}
  For any $\lambda \in \hat{\mX} \subset X$, the character of the projective module 
$P((l-1)\rho + \lambda)$ equals to 
$ch St \cdot f(\lambda)$ in ${\mathbb C}[P]^W$. 
\end{Prop} 
\begin{proof}
By \prop{proj_ch}, the character of $P((l-1)\rho +\lambda)$ in ${\mathbb C}[P]^W$
 can be written as
$$ch P((l-1)\rho +\lambda)= a_\lambda ch St \cdot f(\lambda) + ch St 
\sum_{\mu < \lambda}
a_\mu f(\mu) $$ 
for some coefficients $a_\mu \in {\mathbb C}$.
As a character of an indecomposable projective module over $\Uq$, it contains only Weyl 
characters with highest weights in $\Wl \cdot \lambda$. This requires $\mu \in W \circ \lambda$,
which cannot be satisfied for $\mu < \lambda, \lambda \in \hat{\mX}$. 
Therefore, we have that $ch P((l-1)\rho +\lambda)$ is proportional to 
$ch St \cdot f(\lambda)$,
and the normalization can be checked by evaluating the corresponding characters at $1$. 
 \end{proof}

For any $\lambda \in \hat{\mathcal X}$ chosen as above, 
let $[ {P}(\lambda)]$ denote the character of the restriction of the 
projective $\Uq$-module $P((l-1)\rho +\lambda)$ to $\Uf$. Then 
$\{[ {P}(\lambda)]\}_{\lambda \in \hat{\mX}}$ form a basis in $\bPr$. 

\begin{Cor} For any $\lambda, \mu \in {\hat{\mathcal X}}$, we have in $\R$ 
$$ {[ {P}(\lambda)]} \cdot {[ {P}(\mu)]} = [  {St} ] \cdot 
\sum_{\nu \in {\hat{\mathcal X}}} n_{\lambda, \mu}^\nu {[ {P}(\nu)]},$$ 
where $n_{\lambda, \mu}^\nu$ are the structure constants of the algebra 
$(Vr)^+$ in the basis of elementary $W$-symmetric functions 
$$ f(\lambda) \cdot f(\mu)  = \sum_{\nu \in \hat{\mathcal X}} n_{\lambda, \mu}^\nu f(\nu).$$ 
\end{Cor} 
\begin{proof}
\prop{P(mu)} implies that we have $\bar{\psi}([ {P}(\lambda)])= 
\bar{\psi}([St])\cdot f(\lambda)$ 
in \\ $({\mathbb C}[P]^{W} \otimes_{{\mathbb C}[lP]^W} {\mathbb C})$ for all $\lambda \in \hat{\mX}$, 
and we know that $\bar{\psi}([St])=ch St$ is given by the Weyl character formula. 
Multiplication gives 
$$ {ch St}^2  f(\lambda) f(\mu) = {ch St}^2 \sum_{\nu \in \hat{\mathcal X}} n_{\lambda, \mu}^\nu f_\nu .$$ 
All characters of projective $\Uf$-modules in one $W\bullet$-orbit are proportional, and $\hat{\mathcal X}$ is the fundamental domain for $W\circ$-action which corresponds to the $W\bullet$-action after multiplication by 
$ch St$. Therefore, for each $\nu \in \hat{\mathcal X}$ there exists a unique projective character $[P(\nu)] \in \R$ such that  
${\bar{\psi}([ P(\nu)])} = ch St \cdot f(\nu)$, and using the algebra isomorphism $\bar{\psi}$,
we obtain the Corollary. 
\end{proof}

To summarize, we note that "modulo $[St]$" the ideal $\bPr$ has the structure of a semisimple commutative algebra of dimension $|\hat{\mathcal X}|$, which is isomorphic to the algebra $\R/{\rm Rad}(\R) \simeq {\rm Head}(\R)$ ($\R$ is 
a Frobenius algebra), 
and admits a natural basis of $W$-symmetric functions $f(\lambda)$, $\lambda \in \hat{\mathcal X}$, with each $f(\lambda)$ corresponding to an individual character of an indecomposable projective module. 

To compare, we have for $Vr$: 
$$ [W(\lambda)] [W(\mu)] = \sum a_{\lambda, \mu}^\nu [W(\nu)], $$ 
where $[W(\lambda)], [W(\mu)] \in \mR$ are characters of Weyl modules of highest 
weights $\lambda, \mu \in X$, and  
$a_{\lambda, \mu}^\nu$ are the structure constants of the algebra 
$(Vr)^-$ as defined in \sec{In}.

Now suppose that $l >h$ is even. Then folowing \cite{AP}, it is possible 
to define just as before the big quantum group $\Uq_{ev}$, the category 
of tilting modules, the fusion category and the Verlinde algebra $Vr_{ev}$. 
An important property of the $Vr_{ev}$ is the existence on it 
of a representation of the modular 
group. Define $ch_qW(\lambda)$ to be a 
functional on $\Uq_{ev}$, equal on an element $x \in \Uq$ to the 
trace of $K_{2\rho}x$
 over the Weyl module $W(\lambda)$. Then the linear span of  
$\{ch_q W(\lambda) \}_{\lambda \in X}$ carries a projective representation of the 
modular group (see e.g. \cite{Lyu2}, \cite{BK}). The action is given by the 
bijective mappings $(F_{ss}, T_{ss})$ on  
$\{ch_q W(\lambda) \}_{\lambda \in X}$ which are defined as follows. 
Let 
$$\mu_{ss} = \frac{1}{\sqrt{\sum_{\lambda \in X}{\rm dim}^2_qW(\lambda)}} 
\sum_{\lambda \in X}{\rm dim}_qW(\lambda) ch_q W(\lambda),$$ 
where ${\rm dim}_qW(\lambda)=
 ch_q W(\lambda)(1)$.  Let $R_{12}$ be the canonical element of the quasitriangular 
algebra $\Uq_{ev}$ (see e.g. \cite{LusB}). For any functional $p$ over $\Uq$ define 
$$ { F}_{ss}(p)= \sum_{(\mu_{ss})} (\mu_{ss})_{(1)}(({\rm id}\otimes p)(R_{21}R_{12}))
 (\mu_{ss})_{(2)}, $$
where $\Delta \mu_{ss} = \sum_{(\mu_{ss})}(\mu_{ss})_{(1)} \otimes (\mu_{ss})_{(2)}$ 
is the comultiplication in the Hopf dual of $\Uq_{ev}$. 
Then $F_{ss}$ has the properties of a Fourier transform, its square is proportional to 
the antipode on the Hopf dual of $\Uq_{ev}$. Let $v$ denote the invertible central element 
in $\Uq_{ev}$ such that $v^2 = u\gamma(u)$, where $u = (\gamma \otimes {\rm id})(R_{21})$ and 
$\gamma$ is 
the antipode of $\Uq_{ev}$. Denote by $T_{ss}:  \{ch_q W(\lambda) \}_{\lambda \in X} \to 
\{ch_q W(\lambda) \}_{\lambda \in X}$ the mapping which sends $p(\cdot)$ to $p(v\cdot)$. 
Then by \cite{BK},  ${ F}_{ss}$ and $T_{ss}$ satisfy the modular identity: 
$(F_{ss} T_{ss})^3 = c F_{ss}^2$, where $c \in {\mathbb C}$ is a constant.  

By \cite{Dr}, the algebra of $q$-characters of finite dimensional modules over $\Uq_{ev}$
is isomorphic to the complexified Grothendieck ring of the category of finite dimensional modules, with the isomorphism $\vartheta: p(\cdot) \to p((K_{-2\rho} \cdot)$ sending $ch_q W(\lambda)$ to $\chi(\lambda)$. Then the maps 
${\mathcal F}_{ss} = \vartheta \circ F_{ss} \circ \vartheta^{-1}$,  
${\mathcal T}_{ss} =  T_{ss} $ define a modular 
group action on $Vr_{ev}$.

Back to the situation when $l >h$ is odd and coprime to the determinant 
of the Cartan matrix of $\Lg$, we observe that the ideal 
$\bPr$ is invariant under the {\it quantum Fourier transform} ${\mathcal F}_q$, 
which is a slight modification of the one introduced in \cite{LM}.  
This easily extends to a projective representation of the modular group. 

Define the algebra $\R_r$ of $q^{-1}$-characters as a subalgebra of $\Ud$ generated by 
$\{ch_{q^{-1}}W(\lambda) \}_{\lambda \in P_l}$, where $ch_{q^{-1}}W(\lambda)$ 
is a functional on $\Uf$, equal on an element $x \in \Uf$ to the trace of $K_{-2\rho}x$ over the module $W(\lambda)$. Then the map $\vartheta: p(\cdot) \to p((K_{-2\rho} \cdot)$  provides an isomorphism 
of commutative algebras $\vartheta : \R \to \R_r$. For $p \in \R_r$ the 
element 
$ (p \otimes {\rm id})((\bar{R}_{21}\bar{R}_{12}))$ is central in $\Uf$ \cite{Dr}. 
Let $\lambda_r$ denote the right integral of $\Ud$, which is unique up to a scalar multiple. Then we have 
\begin{Prop} \label{prop:fourier}
Let ${F}_q$ denote the  map 
$$ {F}_q(p)= \sum_{(\lambda_r)} (\lambda_r)_{(1)}((\bar{\gamma})(p \otimes {\rm id})(\bar{R}_{21}\bar{R}_{12})) (\lambda_r)_{(2)}, $$
 where $\bar{\gamma}$ is the antipode in $\Uf$.  Then 
${ F}_q$ maps ${\rm Ann}({\rm Rad}\; \R_r)$ bijectively to itself. 
 \end{Prop}
This is equivalent to Theorem 5.2 in \cite{La}. 
By \cite{Lyu1}, $\Uf$ contains the ribbon element $\bar{v}$, which belongs to the center $\Z \subset \Uf$. Define the mapping $T_q : \R_r \to \R_r$ by the action of the ribbon element $\bar{v} \in \Uf$,
$p(\cdot) \to p(\bar{v} \cdot)$. By \cite{La}, there exists a homomorphism $\phi : \Z \to \R_r$ of $\Z$-modules such that $\phi^{-1}({\rm Ann}({\rm Rad}\; \R_r)) \subset {\rm Ann}({\rm Rad}\; \Z)$. Therefore, any central element acts in ${\rm Ann}({\rm Rad}\; \R_r)$ as multiplication by a scalar. Then $T_q$ preserves  ${\rm Ann}({\rm Rad}\; \R_r)$, and is nondegenerate since it is nondegenerate on $\R_r$. Then by \cite{Lyu1},
 the mappings $F_q, T_q$ satisfy 
the modular identities, and therefore  
${\mathcal F}_q \equiv \vartheta^{-1} \circ F_q \circ \vartheta $ together with 
${\mathcal T}_q \equiv T_q $ defines a modular 
group action on $\bPr$. 

\begin{Rem} 
 We used the algebra of $q^{-1}$-characters here to make it consistent with 
\cite{La}; starting with the quantum Fourier transform defined on the span of  
$q$-characters would lead to the same mapping ${\mathcal F}_q$ on $\bPr$. 
\end{Rem}
 
Below we list the properties of $Vr$ and $\bPr$. 

\hspace{-1cm} \begin{tabular}{l|l}  \label{table}
$Vr$ & $\bPr$ \\
\hline
    &       \\
$Vr \simeq (Vr)^- \simeq {\mathbb C}[P]^W/J$, where  & $\bPr \simeq (ch St)(Vr)^+ \simeq (ch St) {\mathbb C}[P]^W/I,$ \\
$J = \langle (\chi(\lambda) + \chi(s \cdot \lambda)\rangle_{\lambda \in P, {\rm refl.}\, s \in \Wl}$ & 
where $I = \langle (f(\lambda) - f(s (\lambda))\rangle_{\lambda \in P, s \in \WlP}$; \\       
  &  \\
$Vr$ has a basis $\{ [W(\lambda)] \}_{\lambda \in {X}}$ & $\bPr$ has a basis $\{[St] \cdot [f(\lambda)]\}_{\lambda \in \hat{\mX}}$\\
  &  \\
${\rm dim}({Vr}) = |X| ,$ the number  & ${\rm dim}(\bPr) = |\hat{\mathcal X}|=|\tmX|,$ the number  \\
of regular blocks in $\mR$ & of all blocks in $\R$ \\
    &   \\

${Vr}$ is the algebra of characters of the & $\bPr$ is the ideal of characters of  \\
fusion category ${\mathcal{F}us}$ over $\Uq$. & projective modules over $\Uf$.  \\
    &     \\
\hline
    &     \\
For even $l$, $Vr_{ev}$ carries a representation & For odd $l$, $\bPr$ carries a representation \\
of the modular group given by $({\mathcal F}_{ss}, {\mathcal T}_{ss})$ & of the modular group given by $({\mathcal F}_q, {\mathcal T}_q)$ \\
  &   \\
\hline
\end{tabular}

\end{document}